\documentclass[12pt]{article}
\usepackage{amssymb}
\usepackage{amsmath}
\usepackage{graphics}
\usepackage{epsfig}
\usepackage{color}

\newtheorem{thm}{Theorem}
\newtheorem{proposition}[thm]{Proposition}

\begin{document}

\title{A time scales approach to coinfection by opportunistic diseases.}

\author{Marcos Marv\'a$^\star$, Ezio Venturino$^{\dagger}$, Rafael Bravo de la Parra$^\star$\\
$^{\star}$Departamento de F\'isica y Matem\'aticas,\\
Universidad de Alcal\'a,\\
28871 Alcal\'a de Henares, Spain\\
\\
$^{\dagger}$Dipartimento di Matematica ``Giuseppe Peano'',\\
Universit\`a di Torino,\\
via Carlo Alberto 10, 10123 Torino, Italy\\
emails: marcos.marva@uah.es, ezio.venturino@unito.it, rafael.bravo@uah.es
}
\maketitle

\subsection*{Abstract}
Traditional biomedical approaches treat diseases in isolation, but the importance of
synergistic disease interactions is now recognized. As a first step
we present and analyze a simple coinfection model for two diseases affecting simultaneously a population.
The host population is affected by 
the \emph{primary disease}, a
long-term infection whose dynamics is described by a SIS model with demography,
which facilitates individuals acquiring a second disease, 
\emph{secondary (or \emph{opportunistic}) disease}.
The secondary disease is instead a short-term infection affecting only the primary-infected individuals. 
Its dynamics is also represented by a SIS model with no demography.
To distinguish between short and long-term infection the complete model is written as a two 
time scales system. The primary disease acts at the slow time scale while the secondary disease does at 
the fast one, allowing a dimension reduction of the system and making its analysis tractable.
We show that an opportunistic disease outbreak might change drastically the outcome 
of the primary epidemic process, although 
it does among the outcomes allowed by the primary disease. We have found situations in which either 
acting on the opportunistic disease transmission or recovery rates 
or controlling the susceptible and infected population size allow to eradicate/promote disease endemicity.

\section{Introduction.}\label{intro}
Coinfection is the simultaneous infection of a host by multiple pathogen species. 
The global incidence of coinfection among humans is huge \cite{Cox01} and supposed 
to be more common than single infection. The interactions between pathogen species 
within their host can have either positive or negative effects on each other. 
The net effect of coinfection on human health is thought to be negative \cite{Griffiths11}.

The case of positive parasite interactions falls into the concept of \emph{syndemic}: 
aggregation of two or more diseases in a population in which there is some level of 
positive biological interaction that exacerbates the negative health effects of any 
or all of the diseases \cite{Merrill09}. From the point of view of prevention and 
treatment of disease it is the opposite case that is important,
sometimes called \emph{counter syndemic}: 
disease interactions that yield a lower whole effect than the sum effects of the individual 
diseases involved. An example of counter syndemic is that of human immunodeficiency virus (HIV) 
transiently suppressed during acute measles infections. A broadly extended syndemic involves 
tuberculosis (TB) and HIV \cite{Kwan11}. The World Health Organization \cite{WHO} reports that 
people living with HIV are around 30 times more likely to develop TB than persons without HIV 
and also that TB is the most common illness occurring among people living with HIV. Other 
syndemics involving infectious diseases have been described in the literature: HIV and malaria 
syndemic \cite{Abu-Raddad06}; the helminthic infections, malaria and HIV/AIDS syndemic \cite{Ivan13}; 
the pertussis, influenza, tuberculosis syndemic \cite{Herring07}; and the HIV and sexually transmitted 
disease (STD) syndemic \cite{Eaton11}.

In this work we deal with a particular, but very common, type of coinfection. 
We consider the interactions of two diseases, the first one of the type called \emph{primary disease} 
and the second one of the \emph{opportunistic disease} type. Only relatively few pathogen species cause 
disease in otherwise healthy individuals \cite{Baron96}. Those few are called primary pathogens. 
The diseases that they cause, primary diseases, are the result of their only activity within a healthy host. 
An opportunistic disease, on the other hand, is characterized \cite{Symmers65} as a serious, usually 
progressive infection by a micro-organism that has limited (or no) pathogenic capacity under ordinary 
circumstances, but which has been able to cause serious disease as a result of the predisposing effect 
of another disease or of its treatment. 

The importance of opportunistic diseases for public health \cite{Griffiths11,Koster81,Zlamy13} is 
underrepresented in the mathematical modelling literature. A reason for that is that models of 
coinfection usually result in large dimensional systems which are difficult to be studied analytically. 
The main aim of this work is to settle a model describing 
the interaction between both, the primary and the secondary diseases. The model that we present 
in this work tries to capture the basic features of a coinfection model using for it the least 
possible number of variables.  The dynamics of the primary disease is represented by means of a SIS model. 
All individuals affected by the primary infection are assumed to be susceptible of being infected by the 
opportunistic disease. As the dynamics of the opportunistic disease is also described in terms of a SIS 
model we only distinguish three types of individuals in the population: individuals with no infection, 
\emph{susceptible}, individuals infected by the primary disease but not by the opportunistic disease, 
\emph{primary infected}, and individuals infected by both diseases, \emph{coinfected}.

Specifically, we  want to know whether or not the coinfection by a secondary disease produces 
epidemiological scenarios not allowed by the primary disease submodel. In the latter case, 
it is of interest to assess if coinfection has any influence on the actual outcome of the model, even if it is only
among those allowed 
by the primary disease submodel. On the other hand, and in any case, we look for identifying 
mechanisms to modulate the epidemiological outcome. \\

A primary disease enabling secondary infections has typically a long illness period. It must produce a 
persistent alteration of the immune response which weakens the body's ability to clear secondary diseases. 
On the other hand, a compromised immune system presents an opportunity that a secondary pathogen must 
rapidly take advantage of. As a simplified approximation of the general case we suppose that the primary 
disease is a long-term infection that evolves slowly compared to the opportunistic disease which has a rapid 
evolution and, thus, can be considered a short-term infection. This difference in the acting speed of both 
infections is reflected in our model in two different issues. Firstly, we assume that demography has an impact 
in the primary disease, due to its slow evolution, whereas it is negligible for the opportunistic disease which 
evolves in a short period of time. Secondly, the system of differential equations, in terms of which we express 
our model, possesses two time scales: the slow one encompassing the demography and the primary disease evolution 
and the fast one associated to the opportunistic disease evolution.

The inclusion of two time scales in the system has the advantage of allowing its reduction. The asymptotic behavior 
of the solutions of the initial three dimensional system can be studied through a planar system. The reduction of 
the system is undertaken with the help of aggregation methods \cite{Auger08a,Auger12,Metal12}. The general aim of 
these methods is studying the relationships between a large class of complex systems, in which many variables are 
involved, and their corresponding reduced or aggregated systems, governed by a few global variables. The idea behind 
the reduction of the system in our model is considering the evolution of the secondary infection as instantaneous 
in relation to that of the primary one. Obviously this is but an approximation which, on the other hand, can be 
precisely treated with the help of the aggregation method. The steady state rapidly, almost instantaneously, 
reached by the opportunistic disease serves to merging in one single variable those variables corresponding to 
primary infected and coinfected individuals. The result is a SIS type model where the effect of the opportunistic 
disease is reflected in its parameters.\\

The model is presented in Section \ref{mod}. In this section the reduction of the system
is also included. 
Section \ref{anrs} is devoted to the analysis of the reduced system. This analysis allows a discussion of the permanence of 
the population as well as of the influence of the final size of the opportunistic disease on the outcome of the primary epidemic. 
This discussion is the content of section \ref{discu}.

\section{The model}\label{mod}

We build up in this section a model of coinfection that describes the interaction between two diseases, 
one of primary type whereas the second one is of opportunistic type. Only the individuals infected by 
the primary disease are susceptible of being infected by the opportunistic disease. Moreover, the 
interaction of both diseases occurs at different time scales, being the evolution of the opportunistic 
disease much faster than that of the primary one. The model is written in terms of a slow-fast ordinary 
differential equations model. 
After building the slow-fast model, the separation of time scales allows us to apply approximate 
aggregation techniques \cite{Auger08a, Auger12} to 
get a smaller dimensional system. For the convenience of readers non familiar with it, 
he reduction procedure is sketched in subsection \ref{aggrmod}. The section
finishes describing which kind of information about the  
slow-fast system can be retrieved from the reduced system.

\subsection{The primary disease sub-model}\label{1submod}

The primary disease dynamics is described by a SIS model with demographic effects.
In a SIS model individuals are divided into susceptible (S) and infected (I). 
The latter return to the susceptible class on recovery because the disease confers 
no immunity against reinfection \cite{Brauer01}. It is appropriate for most diseases 
transmitted by bacterial or helminth agents, and most sexually transmitted diseases. 
Concerning transmission, there are two extreme traditional forms \cite{Kang14} the 
density-dependent transmission (DDT) and the frequency-dependent transmission (FDT). 
In DDT the rate of contact between susceptible and infected individuals increases with 
host density while in FDT this rate of contact is independent of host density. The fact 
that the primary disease acts together with demography at the same time scale leads us 
to assume it does with density-dependent transmission. On the other hand, in the case of 
the opportunistic disease which turns out to evolve at a faster time scale we consider 
that it does with frequency-dependent transmission \cite{McCallum01,Begon02}.

We denote by $\gamma$ the recovery rate and by $\beta$ the constant transmission rate. 
The parameter $\mu$ describes the additional disease-induced mortality.

We consider demographic effects with only horizontal transmission of the disease. 
In many mathematical models, from a demographic point of view, the differences between 
susceptible and infected individuals are reduced to an additional disease-related death 
rate or disease induced reduction in fecundity \cite{Brauer08}. However, there are experimental 
evidences of the influence of disease on host competitive abilities \cite{Bedhomme05} 
which have already been introduced in eco-epidemiological models \cite{Sieber14}. We adopt 
this last approach. The intrinsic per-capita fertility rate of uninfected individuals is 
given by $r$. The reduction on intrinsic per-capita fertility rate of infected individuals 
is represented by the parameter $a\in (0,1)$. The natural death rate is denoted $m$. 
The effects of intra-specific competition reducing population growth are introduced in 
the model by means of parameters $c_{SS}$, $c_{SI}$, $c_{IS}$ and $c_{II}$. 
To be precise, the parameters $c_{SS}$ and $c_{II}$ represent intra-class competition 
between susceptible and infected individuals, respectively, whereas the parameters 
$c_{SI}$ and $c_{IS}$ introduce the inter-class impact of infected on susceptible individuals
and of susceptible on infected individuals, respectively.

The primary disease submodel is given by the equations
\begin{equation}\label{slow}
\left\{\begin{array}{ccccc}
     \dfrac{dS}{dt} & = & rS+arI-mS & -\left(c_{SS} S+c_{SI}I\right)S & -\beta SI + \gamma I,  \\ \rule{0ex}{5ex}
     	   \dfrac{d I}{dt} & = & \underbrace{ \rule{5ex}{0ex} -m I \rule{5ex}{0ex} }_{\text{\tiny   density independent growth}} & \underbrace{-\left(c_{IS}S+c_{II}I\right)I}_{\text{\tiny   competition}} & \underbrace{+\beta SI - \gamma I -\mu I}_{\text{\tiny   transmission, recovery and disease mortality}}
    \end{array}\right.
\end{equation}
As mentioned in the introduction, to our knowledge, the primary disease submodel (\ref{slow}) has not been previously analyzed. 
However, we postpone its analysis until section \ref{anrs}, once we have described the full model and the aforementioned reduction process.
\subsection{The opportunistic disease sub-model}\label{2submod}
The opportunistic disease spreads only through the individuals infected by the primary disease.
We consider that the opportunistic disease dynamics is also described by a SIS model. 
Individuals infected by the primary disease are further classified into those not infected 
by the opportunistic disease (U), primary infected, and those infected by both diseases (V), coinfected.

The fast evolution of the opportunistic disease, compared to primary disease and demography, 
suggests not including demographic effects and choosing the frequency-dependent transmission 
form. Let $\lambda$ and $\delta$ be, respectively, the constant transmission and recovery rates.

The opportunistic disease submodel is represented by the equations
\begin{equation}\label{fast}
\left\{\begin{array}{cccccc}
	   \dfrac{d U}{d \tau} & = & - & \lambda \dfrac{UV}{U+V} & + & \delta V,\\ \rule{0ex}{5ex}
     	   \dfrac{d V}{d \tau}  & = &  &  \lambda \dfrac{UV}{U+V} & - &  \delta V.
            \end{array}
\right.
\end{equation}
We use $\tau$ to denote the time variable for the fast time scale. It is related to variable time $t$  
in system (\ref{slow}) as $t=\varepsilon \tau$ where $\varepsilon$ is a small positive constant 
representing the ratio between time scales.

\subsection{The full two time-scales model}\label{fullmod}

Finally, we construct  the model encompassing both diseases. It has the form of a system with 
three state variables: susceptible $S$, primary infected $U$ and coinfected $V$ individuals. 
It is a system with two time scales that is expressed  in terms of the fast time variable $\tau$. 
The terms associated to the slow time scale, demography and primary disease dynamics, appear multiplied 
by $\varepsilon$ in (\ref{complete}). The fast part of system (\ref{complete}), the opportunistic disease
dynamics, coincides with system (\ref{fast}).

In the slow part of  system (\ref{complete}) we have to define different rates for primary infected 
and coinfected individuals.
We denote $\beta_U$ and $\beta_V$ the constant primary disease transmission rates due to primary 
infected and coinfected individuals, respectively.
We assume that there is no direct connection between 
the susceptible and coinfected stages. A susceptible 
individual must first acquire the primary disease and 
later be infected by the opportunistic one. On the other 
hand, a coinfected individual must first recover from the 
opportunistic disease and then, being just primary infected, can also recover from the primary one. 
The primary disease recovery rate is still denoted $\gamma$.   Parameters $\mu_U$ and $\mu_V$ describe 
the additional  primary disease-induced mortality in primary infected and coinfected individuals, respectively.

Concerning the part of demography, we keep the same intrinsic per-capita fertility rate
of uninfected individuals $r$ and the natural death rate of the population $m$
as in system (\ref{slow}). 
We include different coefficients of reduction on intrinsic per-capita fertility rate for primary 
infected and coinfected individuals: $a_U$ and $a_V$. We assume them to verify $0<a_V<a_U<1$ 
supposing that coinfected individuals participate in reproduction though at a smaller rate. 
To distinguish the effects of intraspecific competition among the three stages we need to 
introduce nine parameters $c_{SS}$, $c_{SU}$, $c_{SV}$, $c_{US}$, $c_{UU}$, $c_{UV}$, $c_{VS}$, 
$c_{VU}$ and $c_{VV}$. They represent the competition, either intra-class o inter-class, between 
the two stages in each of the nine different interaction pairs.

The complete two time scales system reads as follows
\begin{equation}\label{complete}
\left\{
\begin{array}{l}
     \dfrac{dS}{d\tau}= \varepsilon\left[\rule{0ex}{2.5ex}rS+a_UrU+a_VrV-mS-(c_{SS}S+c_{SU}U+c_{SV}V)S
	      -\beta_U SU-\beta_V SV  + \gamma U\right],\\  \rule{0ex}{5ex}
     \dfrac{d U}{d\tau}= -\dfrac{\lambda UV}{U+V} + \delta V \\ \rule{0ex}{3ex}
  \rule{14ex}{0ex} + \varepsilon\left[\rule{0ex}{2.5ex}-mU-(c_{US}S+c_{UU}U+c_{UV}V)U+\beta_U SU+\beta_V SV  - \gamma U-\mu_U U\right],\\  \rule{0ex}{5ex}
	   \dfrac{d V}{d\tau}= \dfrac{\lambda UV}{U+V} - \delta V+
	    \varepsilon\left[\rule{0ex}{2.5ex}-mV-(c_{VS}S+c_{VU}U+c_{VV}V)V-\mu_V V\right].\\
            \end{array}
\right.
          \end{equation}

\subsection{Reduction of the model}\label{aggrmod}

In this section we take advantage of the two time scales to reduce the dimension of the
complete  system (\ref{complete}). In the next section, as a consequence of this reduction, 
we perform the analysis of the model by means of a planar system. The reduction follows 
the technique called \emph{approximate aggregation method} \cite{Auger08a,Auger12}. The first 
step is writing the system in the so-called slow-fast form. This is easily done in system 
(\ref{complete}) using the change of variables $(S,U,V)\mapsto(S,I,V)$ where $I=U+V$ 
represents all the infected individual, both primary infected and coinfected.
\begin{equation}\label{slow-fast}
\left\{
\begin{array}{l}
     \dfrac{dS}{d\tau}= \varepsilon\left[\rule{0ex}{2.5ex}rS+a_Ur(I-V)+a_VrV-mS-(c_{SS}S+c_{SU}(I-V)+c_{SV}V)S\right. \\\rule{0ex}{3ex}
  \rule{47ex}{0ex}    \left.\rule{0ex}{2.5ex}   -\beta_U S(I-V)-\beta_V SV  + \gamma (I-V)\right],\\  \rule{0ex}{5ex}
\dfrac{d I}{d\tau}=\varepsilon\left[\rule{0ex}{2.5ex}-mI-(c_{US}S+c_{UU}(I-V)+c_{UV}V)(I-V)+\beta_U S(I-V)+\beta_V SV \right. \\\rule{0ex}{3ex}
 \rule{16ex}{0ex}    \left.\rule{0ex}{2.5ex}- \gamma (I-V)-\mu_U (I-V)-(c_{VS}S+c_{VU}(I-V)+c_{VV}V)V-\mu_V V \right],\\  \rule{0ex}{5ex}
	   \dfrac{d V}{d\tau}= \dfrac{\lambda (I-V)V}{I} - \delta V+
	    \varepsilon\left[\rule{0ex}{2.5ex}-mV-(c_{VS}S+c_{VU}(I-V)+c_{VV}V)V-\mu_V V\right].\\
            \end{array}
\right.
          \end{equation}
The key point of the new form of system (\ref{complete}) is making visible that variables $S$ and $I$ 
are slow (the right hand side terms of their equations have $\varepsilon$ as a factor) in the sense that they almost do not change at the fast time scale. The fast dynamics is concentrated in the first terms without $\varepsilon$ in the equation for $V$. The approximation that the aggregation method proposes consist in separating both dynamics. Firstly, the non slow variables are calculated in terms of the slow ones by assuming that they are the equilibria (called \emph{fast equilibria}) determining the long term behaviour of the fast dynamics. Secondly, these obtained values of the non slow variables are substituted into the equations of the slow ones yielding a reduced system for the latter. In this reduced or aggregated system the fast dynamics is summarized 
in its parameters. In the particular case of system (\ref{slow-fast}), the only non slow variable is $V$ and the fast dynamics reduces to the equation
\begin{equation}\label{fasteq}
	   \dfrac{d V}{d\tau}= \dfrac{\lambda(I-V)V}{I} - \delta V.
\end{equation}
Assuming the slow variable $I$ to be constant, the analysis of equation (\ref{fasteq}) gives that for positive values of $V(0)$
\begin{equation}\label{fasteq1}
    \lim_{\tau\to\infty}V(\tau)=\nu^*I=\left\{\begin{array}{ll}
                            0 & \text{if }\delta\geq \lambda\\
				   (1-\delta/\lambda)I & \text{if }\delta<\lambda
                           \end{array}\right.
\end{equation}
which corresponds to the results of a classical SIS model without demography and frequency-dependent transmission \cite{Brauer01}. If the recovery rate is larger than the transmission rate, the disease disappears since the number of coinfected individuals tends rapidly to zero. On the other hand, if the recovery rate is smaller than the transmission rate, the disease becomes endemic with a stable fraction $\nu^*=1-\delta/\lambda$ of the infected population $I$ remaining coinfected.

The fast equilibria  $V=\nu^*I$ found in (\ref{fasteq}) are the values to be substituted into the equations for the slow variables $S$ and $I$ to obtain the following reduced system
\begin{equation}\label{redsys}
\left\{
\begin{array}{l}
     \dfrac{dS}{dt}=rS+\bar{a}rI-mS -\left(c_{SS} S+\bar{c}_{SI}I\right)S -\bar{\beta} SI + \bar{\gamma} I,  \\ \rule{0ex}{5ex}
     	   \dfrac{d I}{dt} = -m I -\left(\bar{c}_{IS}S+\bar{c}_{II}I\right)I+\bar{\beta} SI - \bar{\gamma} I -\bar{\mu} I,
    \end{array}
\right.
\end{equation}
which has the same form as the primary disease sub-model (\ref{slow}). In its parameters
the effect of fast dynamics, the opportunistic disease, is implicit through $\nu^*$:
\begin{equation}\label{coeffs}
 \begin{array}{l}
\bar{a}=(1-\nu^*)a_U+\nu^*a_V,\ \bar{c}_{SI}=(1-\nu^*)c_{SU}+\nu^*c_{SV},\ \bar{c}_{IS}=(1-\nu^*)c_{US}+\nu^*c_{VS},  \\ \rule{0ex}{3ex}
\bar{c}_{II}=(1-\nu^*)^2c_{UU}+(1-\nu^*)\nu^*c_{UV}+\nu^*(1-\nu^*)c_{VU}+(\nu^*)^2c_{VV}, \\ \rule{0ex}{3ex}
 \bar{\beta}=(1-\nu^*)\beta_U+\nu^*\beta_V,\ \bar{\gamma}=(1-\nu^*)\gamma,\ \bar{\mu}=(1-\nu^*)\mu_U+\nu^*\mu_V
 \end{array}
\end{equation}
The reduced system (\ref{redsys}) is useful to analyze the asymptotic behaviour of the solutions 
of the complete system (\ref{complete}) \cite{Auger08a,Auger12}. In particular, the 
existence of a hyperbolic asymptotically stable equilibria $(S^*,I^*)$ of system (\ref{redsys}) 
ensures the existence, for $\varepsilon$ small enough, of an equilibria of system (\ref{complete}) 
with the same characteristics and a form very close to $(S^*,(1-\nu^*)I^*,\nu^*I^*)$. 
In the next section we carry out the analysis of the stability of equilibria of system (\ref{redsys}) 
obtaining thus the corresponding results for the complete model (\ref{complete}).\newline

Note that the reduced system (\ref{redsys}) and the 
primary disease submodel (\ref{slow}) are the same, the only difference being the values of the 
respective coefficients. Indeed, when $\delta>\lambda$ the opportunistic disease cannot invade the population and, in this case, the coefficients (\ref{coeffs}) of systems (\ref{redsys}) 
and (\ref{slow})  are exactly the same.

\section{Analysis of the reduced system}\label{anrs}

We proceed in this section to analyze the reduced system (\ref{redsys}).

We first note that $E_0^*=(0,0)$ is an equilibrium point, the positive $S$ semi-axis, $\{(S,0):S>0\}$, is invariant and on the positive $I$ semi-axis , $\{(0,I):I>0\}$, the vector field associated to system (\ref{redsys}) 
points to the interior of the positive quadrant. We then have that the closed positive quadrant $\mathbb{R}^2_+=\{(S,I):S\geq 0,I\geq 0\}$ is positively invariant.

In the next result we prove that, as expected, if the susceptible fertility rate $r$ is not strictly larger than the natural death rate $m$  the population gets extinct.

\begin{proposition}\label{prop1} If $r\leq m$ then any solution $(S(t),I(t))$ of system (\ref{redsys}) with non-negative initial conditions $(S(0),I(0))$ tends to $E_0^*$.
\end{proposition}

\textbf{Proof}: Let us call $W=S+I$.  Summing up both equations in system (\ref{redsys}) 
and now choosing $c=\min\{c_{SS} ,\bar{c}_{SI},\bar{c}_{IS},\bar{c}_{II}\}$ we have $\dfrac{dW}{dt}\leq  -cW^2$
that, by integration, yields
$$
0\leq W(t)\leq \frac{W(0)}{1+W(0)ct}\underset{t\rightarrow \infty}{\longrightarrow}0
$$
since $\mathbb{R}^2_+$ is positively invariant. 
$\blacksquare$\\\newline

Henceforth we assume that $r>m$. This assumption prevents the population from extinction. The linearization of system (\ref{redsys}) at the equilibrium $E_0^*$ has the matrix
$$
\left(
    \begin{array}{cc}
      r-m & \bar{a}r+\bar{\gamma} \\ \rule{0ex}{3ex}
      0 & -(m+\bar{\gamma}+\bar{\mu})
    \end{array}
  \right)
$$
with one positive and one negative eigenvalues. The unstable manifold of $E_0^*$, associated to $r-m$, is included in the $S$ axis, while the stable manifold, 
associated to $-(m+\bar{\gamma}+\bar{\mu})$, is tangent at $\mathbf{0}$ to the eigenvector $\left(\bar{a}r+\bar{\gamma}, -(r+\bar{\gamma}+\bar{\mu})\right)$ and lies completely
outside the interior of the positive quadrant.

Assuming $r>m$ the only non-negative solution tending to $\mathbf{0}$ is $E_0^*$ itself. We prove next that all non-negative solutions are forward bounded.
\begin{proposition}\label{prop2} Let $r>m$. If $(S(t),I(t))$ is any solution of system (\ref{redsys}) with non-negative initial conditions $(S(0),I(0))$ then it is bounded on $[0,\infty)$.
\end{proposition}
\textbf{Proof}: Calling $W=S+I$ 
%we have
%$$
%\dfrac{dW}{dt}=(r-m)S+(\bar{a}r-m)I -\left(c_{SS} S+\bar{c}_{SI}I\right)S -\left(\bar{c}_{IS}S+\bar{c}_{II}I\right)I-\bar{\mu} I.
%$$ 
and letting 
%Let 
$c=\min\{c_{SS} ,\bar{c}_{SI},\bar{c}_{IS},\bar{c}_{II}\}$ we have
$$
\dfrac{dW}{dt}+(r-m)W\leq 2(r-m)W-cW^2
$$
Function $g(W)=2(r-m)W-cW^2$ attains its maximum on $[0,\infty)$ at $W=(r-m)/c$, 
%$g((r-m)/c)=(r-m)^2/c$.
so that 
%Now we have
$$
\dfrac{dW}{dt}+(r-m)W\leq \dfrac{(r-m)^2}{c}.$$ 
Multiplying both sides of the previous inequality by $e^{(r-m)t}$ and rearranging terms 
%$$e^{(r-m)t}\dfrac{dW}{dt}+(r-m)e^{(r-m)t}W\leq e^{(r-m)t}\dfrac{(r-m)^2}{c}
% $$
%which is equivalent to 
yields 
 $$
 \dfrac{d}{dt}(e^{(r-m)t}W)\leq e^{(r-m)t}\dfrac{(r-m)^2}{c},
$$
that implies, integrating on $[0,t]$,
$$
e^{(r-m)t}W(t)-W(0)\leq \dfrac{(r-m)^2}{c(r-m)}(e^{(r-m)t}-1).$$ 
Rearranging terms in the previous expression leads to 
$$W(t)\leq W(0)e^{-(r-m)t}+\dfrac{r-m}{c}(1-e^{-(r-m)t})
$$
and, finally, we have that $W(t)\leq \max\{W(0),(r-m)/c\}$ for every $t\in [0,\infty)$.
$\blacksquare$\\\newline

In addition to the trivial equilibrium $E_0^*$, the system (\ref{redsys}) possesses a disease-free equilibrium $E_1^*=(S^*_1,0)$, where
$$S^*_1=\frac{r-m}{c_{SS}}$$
that represents the stable size of the population in case of no infection. The growth of the population in the absence of infections is logistic and
$S^*_1$ represents its carrying capacity.
\begin{proposition}\label{prop3} Let $r>m$. The equilibrium point $E_1^*=(S^*_1,0)$ of system (\ref{redsys}) verifies:
\begin{enumerate}
  \item If $S^*_1(\bar{\beta}-\bar{c}_{IS})>m+\bar{\gamma}+\bar{\mu}$ then $E_1^*$ is a saddle point of which stable manifold coincides with the positive $S$ semi-axis.
  \item If $S^*_1(\bar{\beta}-\bar{c}_{IS})<m+\bar{\gamma}+\bar{\mu}$ then $E_1^*$ is locally asymptotically stable.
  \item If $\bar{\beta}-\bar{c}_{IS}\leq 0$ then the basin of attraction of $E_1^*$ includes $\mathbb{R}^2_+-\{\mathbf{0}\}$.
\end{enumerate}
\end{proposition}
\textbf{Proof}:
To prove the two first items it suffices to calculate the matrix of the linearization of system (\ref{redsys}) at $E_1^*$
$$
\left(
    \begin{array}{cc}
      -(r-m) & \bar{a}r+\bar{\gamma}-S^*_1(\bar{c}_{SI}+\bar{\beta}) \\\rule{0ex}{3ex}
      0 & S^*_1(\bar{\beta}-\bar{c}_{IS})-(m+\bar{\gamma}+\bar{\mu})
    \end{array}
  \right)
$$
One of the eigenvalues, $-(r-m)$, is negative while stability depends on the other one, $S^*_1(\bar{\beta}-\bar{c}_{IS})-(m+\bar{\gamma}+\bar{\mu})$, being positive or negative:
$E_1^*$ is a saddle or (locally) asymptotically stable, respectively.

To prove the last assertion we first note that there exists no  interior equilibria because the right-hand side of the $I$ equation is always negative for positive $S$ and $I$. Now the Poincar\'e-Bendixson theorem implies that there is no closed orbit in the interior of the positive quadrant and therefore that all positive solutions, that are forward bounded, must tend to the unique non-negative equilibrium point, $E_1^*$.
$\blacksquare$\\\newline
Up to now we have obtained the condition of non-extinction of the population, $r>m$, and a sufficient condition, $\bar{\beta}-\bar{c}_{IS}\leq 0$, for a disease-free scenario in the long term. This last condition says that if the competition coefficient $\bar{c}_{IS}$, representing the impact of susceptible on infected individuals, is larger that the transmission rate $\bar{\beta}$ then the infection disappears independently of the initial conditions. More significant cases exist when the simple competitive pressure of susceptible on infected individuals is not enough to compensate transmission.

From now on we are also assuming that $\bar{\beta}-\bar{c}_{IS}>0$. In this case the conditions of local stability of the equilibrium $E_1^*$ can be expressed in terms of the parameter
\begin{equation}\label{apar}
    \bar{A}=\frac{m+\bar{\gamma}+\bar{\mu}}{\bar{\beta}-\bar{c}_{IS}}.
\end{equation}
Thus, Proposition \ref{prop3} can be restated as follows: for $r>m$ and $\bar \beta>\bar c_{IS}$, if $S^*_1>\bar{A}$ or $S^*_1<\bar{A}$ 
then  the equilibrium  $E_1^*$ is a saddle point 
or locally asymptotically stable, respectively. The parameter $\bar{A}$ represents a threshold population size allowing or 
not the increase of the infection when it is rare. If the susceptible population is close to its carrying capacity, $S_1^*$, 
a few infected individuals are able to spread the disease if the number of susceptible individuals is large enough, $S^*_1>\bar{A}$. 
On the other hand, the infection disappears if the susceptible population is under the threshold $\bar{A}$.

In the next results we search for conditions ensuring the endemicity of the infection. To express them in a simpler form we define another parameter
 \begin{equation}\label{bpar}
    \bar{B}=\frac{\bar{a}r+\bar{\gamma}}{\bar{c}_{SI}+\bar{\beta}},
\end{equation}
that can be interpreted through the terms depending on $I$ in the first equation of system (\ref{redsys}). This equation can be written in the following form
$$
\dfrac{dS}{dt}=(r-m)S-c_{SS} S^2+\left(\bar{a}r+ \bar{\gamma}-(\bar{c}_{SI}+\bar{\beta})S\right)I
$$
where we note that depending on whether $\bar{a}r+ \bar{\gamma}-(\bar{c}_{SI}+\bar{\beta})S(t)$ is positive or negative, the existence of infected individuals makes the susceptible growth rate increase or decrease, respectively. The size of the susceptible population determines if the infection has a positive or a negative effect on its growing. If $S(t)<\bar{B}$, the more infected individuals the larger the susceptible population growth rate, while $S(t)>\bar{B}$ yields a larger decrease of the susceptible population growth rate whenever there is a larger infected population.

Using parameters $S_1^*$, $\bar{A}$ and $\bar{B}$ the system (\ref{redsys}) can be expressed as follows:
\begin{equation}\label{redsys2}
\left\{
\begin{array}{l}
     \dfrac{dS}{dt}=c_{SS}S(S_1^*-S)+(\bar{c}_{SI}+\bar{\beta})(\bar{B}-S) I,  \\ \rule{0ex}{5ex}
     	   \dfrac{d I}{dt} = (\bar{\beta}-\bar{c}_{IS})(S-\bar{A})I-\bar{c}_{II}I^2.
    \end{array}
\right.
\end{equation}

The equation of the $S$-nullcline of the system (\ref{redsys2}), for $S_1^*\neq \bar{B}$, is
$$
I=\Phi(S)=\frac{c_{SS}}{\bar{c}_{SI}+\bar{\beta}}\cdot \frac{S(S_1^*-S)}{S-\bar{B}}.
$$
We are interested in the part included in the positive quadrant. 
This is for $S_1^*< \bar{B}$  an increasing branch going from the point $E_1^*=(S_1^*,0)$ 
to the asymptote $S=\bar{B}$ (see Figure 1, left panel) and for $S_1^*> \bar{B}$  a decreasing 
branch going from the asymptote $S=\bar{B}$ to the point $E_1^*$ (see Figure 1, right panel). 
In the case $S_1^*= \bar{B}$ the $S$-nullcline in the positive quadrant reduces to 
the line $S=\bar{B}$ (see Figure 1, center panel).

% The $I$-nullcline of the system (\ref{redsys2}) in the open positive quadrant is the line
$$
I=\Psi(S)=\frac{\bar{\beta}-\bar{c}_{IS}}{\bar{c}_{II}}(S-\bar{A})
$$
It is immediate to prove that if $\bar{A}\geq \max\{S_1^*,\bar{B}\}$ there is no  interior equilibria of the system (see panels in Figure 1) and thus, applying again the Poincar\'e-Bendixson theorem, 
we get that all positive solutions tend to $E_1^*$. 
We  gather these results in the next proposition.
\begin{figure}[h!]
 \begin{center}
  \includegraphics[width=13cm]{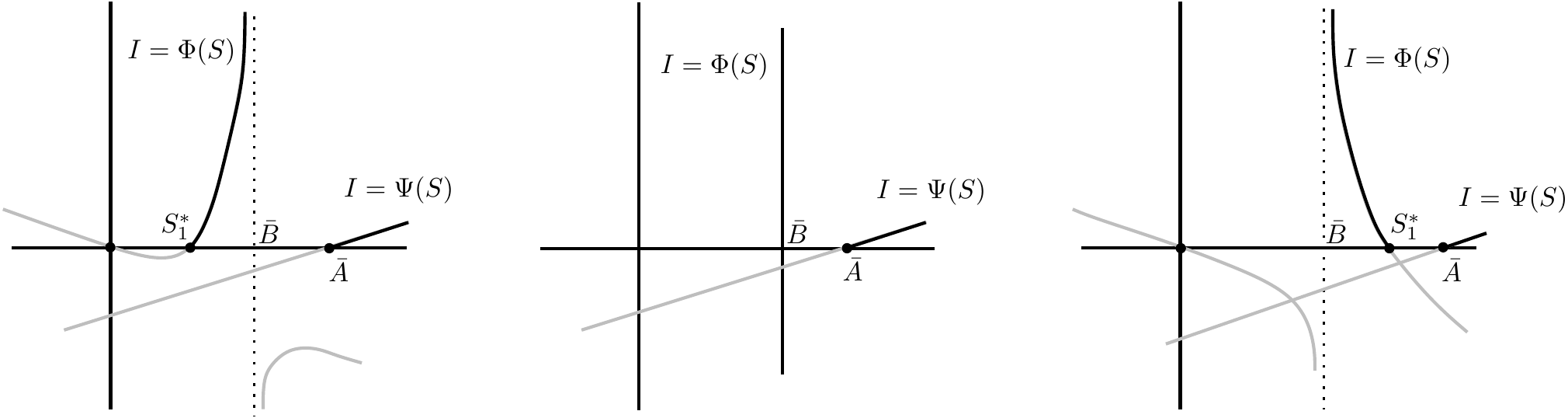}
\label{phases}
 \end{center}
\caption{Possible profiles of the S-nullcline $I=\Phi(S)$ 
and the I-nullcline $I=\Psi(S)$ and the S-nullcline $I=\Phi(S)$  (in black, their intersection with the positive cone).}
\end{figure}

\begin{proposition}\label{prop4} Let $r>m$ and $\bar{\beta}-\bar{c}_{IS}>0$. If  $\bar{A}\geq \max\{S_1^*,\bar{B}\}$ then the system (\ref{redsys}) 
possesses a unique non-negative equilibrium point $E^*_1=(S^*_1,0)$ that is 
asymptotically stable and attracts every positive solution.
\end{proposition}

Condition $\bar{A}\geq \max\{S_1^*,\bar{B}\}$ tell us, on the one hand, that the infection cannot invade due to  $\bar{A}\geq S_1^*$ and, on the other hand, that infected individuals cannot help in attaining the invasion threshold because $\bar{A}\geq \bar{B}$. The consequence is that infection disappears.

There are two situations for the infection to become endemic. The first one is allowing invasion, that is $S_1^* > \bar{A}$, that is treated in 
Proposition \ref{prop5}. The second one does not allow infection invasion for a low number of infected individuals, $S_1^* < \bar{A}$, but larger numbers of infected individuals might help the susceptible population growing, $\bar{A}\geq \bar{B}$, so as to maintain this latter over the invasion threshold. 
In Proposition \ref{prop6} are detailed sufficient conditions to meet this second situation.

\begin{proposition}\label{prop5} Let $r>m$ and $\bar{\beta}-\bar{c}_{IS}>0$. If  $S^*_1>\bar{A}$ then the system (\ref{redsys}) possesses a unique 
 interior equilibrium point $E^*_+=(S^*_+,I^*_+)$ that is locally asymptotically stable. If, in addition, $ S_1^*< \bar B$, then $E^*_+$ attracts every positive  solution.
\end{proposition}
\textbf{Proof}: The assumptions on parameters yield the existence of a unique  
interior equilibrium (see Figure 2). 
The asymptotic stability that follows can be proved by linearization.

Note that the condition $S_1^* > \bar{A}$ might not ensure that all positive solutions tends to 
the interior equilibrium. Due to the Poincar\'e-Bendixson theorem, it might happen that some of 
these solutions tend to a limit cycle included in the open positive quadrant surrounding the equilibrium. 

Condition $S_1^*<\bar B$ excludes the existence of any limit cycle because the region $\mathcal A=\left\{(S,I);\,\Phi(S)\leq S\leq \Psi(S),\,\,0\leq I\right\}$ 
is an invariant (''trapping'') subregion $\mathcal A\subset\mathbb{R}_+^2$ such that 
$E_+^*\in\partial \mathcal A$, the boundary of $\mathcal A$, and 
$\partial \mathcal A\cap\partial \mathbb{R}_+^2\neq \emptyset$ (see Figure 2). That is, any orbit surrounding $E_+^*$ must enter $\mathcal A$, but cannot leave from there. 
\begin{figure}[h]
\begin{center}
  \includegraphics[width=.3\textwidth]{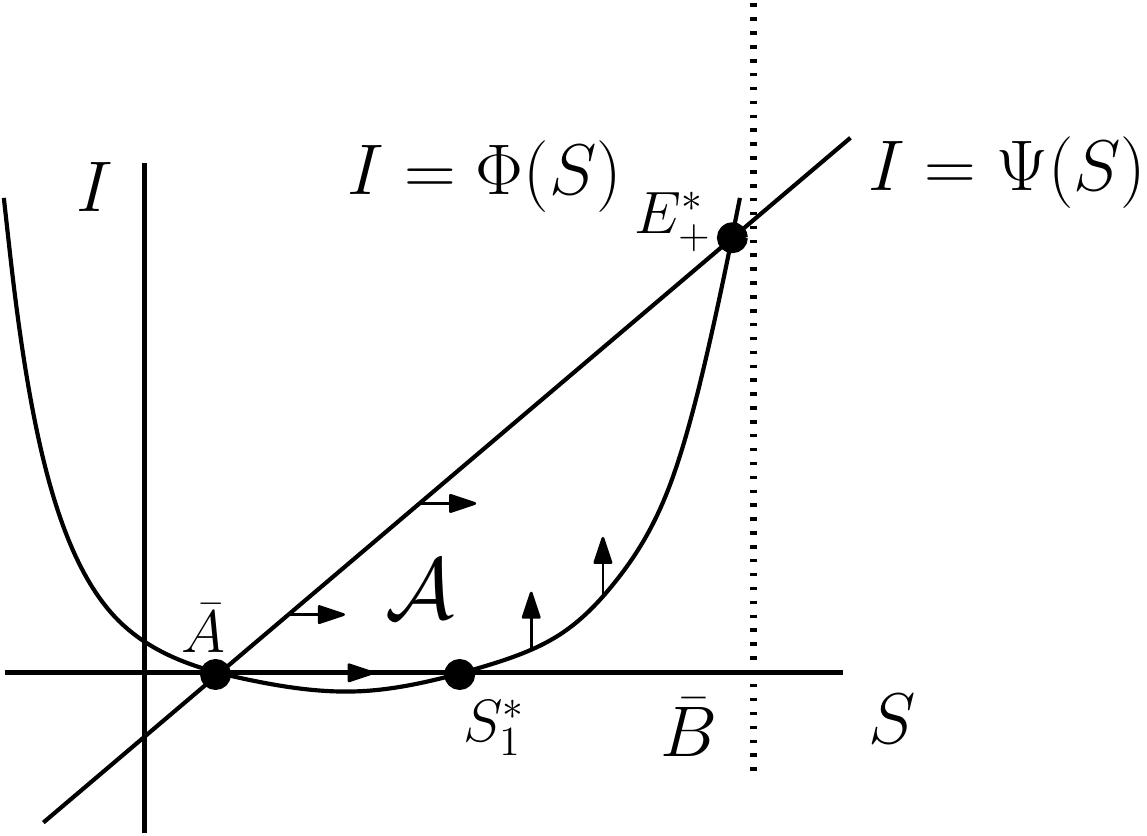}
\end{center}
\label{invariant1}\caption{The invariant region $\mathcal A=\left\{(S,I);\,\Phi(S)\leq S\leq \Psi(S),\,\,0\leq I\right\}$ mentioned in the proof
 of proposition 5.}
\end{figure}
$\blacksquare$\\\newline

In any case, what condition $S_1^* > \bar{A}$ ensures is the endemicity of the infection. 
In the next Proposition \ref{prop6} we state conditions leading the population to the disease-free state or toward conditional endemicity related with a bistable 
scenario.
By conditional endemicity we mean that 
the outcome of the model can be either disease-free (see left panel in Figure 3) or with an endemic disease depending on the initial amount of susceptible and infected individuals; 
see the right panel in Figure 3. 
\begin{figure}[h]
\begin{center}
   \includegraphics[width=.25\textwidth]{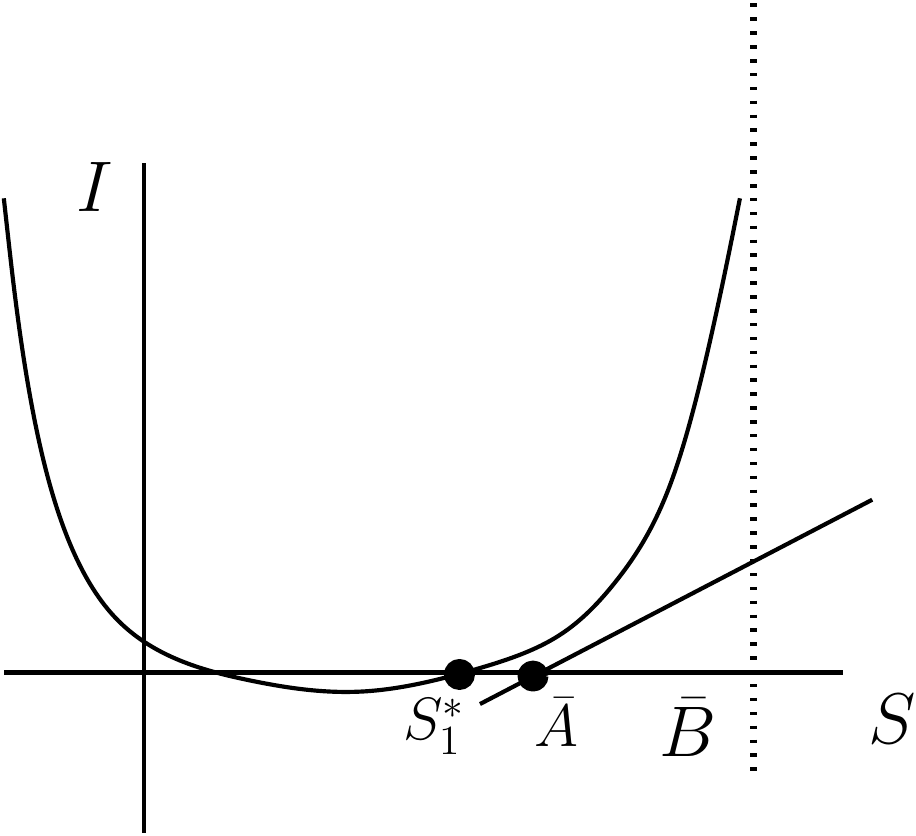}\qquad\qquad
   \includegraphics[width=.25\textwidth]{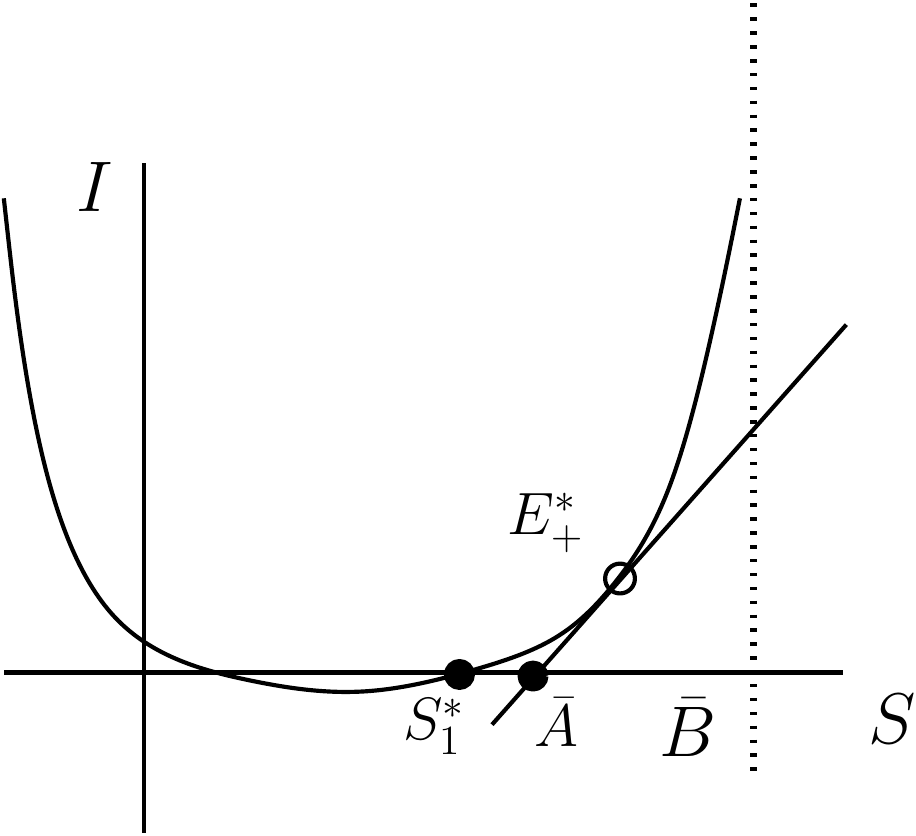}\qquad\qquad
   \includegraphics[width=.25\textwidth]{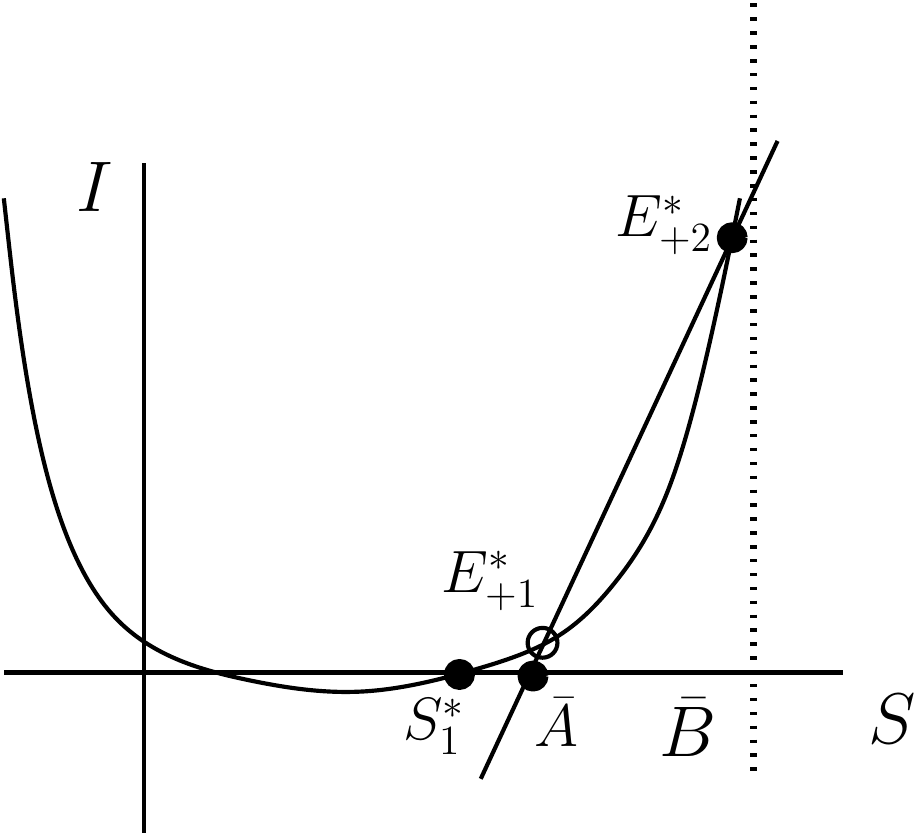}
\end{center}
\label{p6}\caption{Related to Proposition 6, possible configurations of the nullclines when 
$S_1^*<\bar A<\bar B$. On the left panel $S_1^*$ is global attractor. 
On the center panel $S_1^*$ is an attractor. 
On the right panel the bistable case: $S_1^*$ and $E_{+2}^*$ are locally asymptotically stable  while $E_{+1}^*$ is a saddle.}
\end{figure}

Indeed, we introduce $\mathcal R$, that appears later in the corresponding proof, and drives the epidemic outcome. 
This quantity depends on the parameters of the model and allows 
to discriminate whether conditional disease endemicity is allowed or cannot occur.

\begin{proposition}\label{prop6} Let $r>m$ and $\bar{\beta}-\bar{c}_{IS}>0$. If  $S^*_1<\bar{A}<\bar{B}$ then the system (\ref{redsys}) possesses the asymptotically stable equilibrium $E^*_1=(S^*_1,0)$.  
Furthermore, let us define
\[\mathcal R =
\dfrac{\bar{c}_{II}(r-m)+(\bar{a}r+\bar{\gamma})(\bar{\beta}-\bar{c}_{IS})
+(\bar{c}_{SI}+\bar{\beta})(\bar{\mu}+\bar{\gamma})}
{2\sqrt{\left(c_{SS}\bar{c}_{II}+(\bar{c}_{SI}+\bar{\beta})(\bar{\beta}-\bar{c}_{IS})\right)(\bar{a}r+\bar{\gamma})(\bar{\mu}+\bar{\gamma})
}}.
\]
We have:
\begin{enumerate}
  \item if $\mathcal R<1$ then there is no interior equilibrium point and the basin of attraction of $E_1^*$ includes $\mathbb{R}_+^2\backslash\left\{\mathbf{0}\right\}$.
  \item If $\mathcal R = 1$ then there is only one interior equilibrium point $E^*_+=(S^*_+,I^*_+)$  that 
    is unstable. The equilibrium $E^*_1=(S^*_1,0)$ attracts every solution with initial values in the interior of $\mathbb{R}_+^2\backslash E^*_+$. 
  \item If $\mathcal R>1$ then there are two interior equilibrium points $E^*_{+1}=(S^*_{+1},I^*_{+1})$ and $E^*_{+2}=(S^*_{+2},I^*_{+2})$, 
  with $S^*_{+1}<S^*_{+2}$ and $I^*_{+1}<I^*_{+2}$. $E^*_{+1}$ is a saddle and $E^*_{+2}$
is locally asymptotically stable.
\end{enumerate}
\end{proposition}
\textbf{Proof}: Note that the asymptotic stability of $E^*_1$ follows directly from Proposition \ref{prop3} since $S^*_1<\bar{A}$. 

Next, we focus on showing the relation between $\mathcal R$ and the existence of equilibrium points of system (\ref{redsys}). 
Equating the nullclines $I=\Phi(S)$ and $I=\Psi(S)$ of system (\ref{redsys}) yields 
$$
\frac{c_{SS}}{\bar{c}_{SI}+\bar{\beta}}\cdot \frac{S(S_1^*-S)}{S-\bar{B}}=\frac{\bar{\beta}-\bar{c}_{IS}}{\bar{c}_{II}}(S-\bar{A}).
$$
Keeping in mind the definition of $\bar A$ and $\bar B$, the previous expression in equivalent to 
\begin{equation}\label{roots}
\begin{array}{rl}
0=&
\left(c_{ss}\bar{c}_{II}+(\bar{c}_{SI}+\bar{\beta})(\bar{\beta}-\bar{c}_{IS})\right)S^2\\
\rule{0ex}{5ex}
& -
\left(\bar{c}_{II}(r-m)+(\bar{a}r+\bar{\gamma})(\bar{\beta}-\bar{c}_{IS})
+(\bar{c}_{SI}+\bar{\beta})(\bar{\mu}+\bar{\gamma})\right)S \\
\rule{0ex}{5ex}&+
(\bar{a}r+\bar{\gamma})(\bar{\mu}+\bar{\gamma})
\end{array}
\end{equation}
Now, direct calculations lead to the fact that $\mathcal{R}$ being
smaller than, equal to or larger than $1$ 
is equivalent to the discriminant of equation \eqref{roots} 
being negative, zero or positive. This yields the number of equilibrium points.

Note that when $\mathcal R<1$ there is no interior equilibrium point and a direct application of the Poincar\'e-Bendixson theorem yields 
statement (1).

Concerning statement (2), direct calculations show that when $\mathcal R=1$ the equilibrium $E^*$ is not hyperbolic, so that we cannot use the linearization criterion. 
Note that the region $\mathcal B = \left\{(S,I);\,\Psi(S)\leq S\leq \Phi(S),\,\,0\leq I\right\}$ is an invariant region 
such that $E_+^*\in\partial \mathcal B$ and $\partial \mathcal B\cap\partial \mathbb{R}_+^2\neq \emptyset$. Indeed, any solution 
with initial values in the interior of $\mathcal B$ converges to $E_0^*$ since $\Psi(S)< S < \Phi(S)$, which implies that $E_+^*$ is unstable. 
The non existence of periodic orbits can be argued as done in the proof of Proposition \ref{prop5}.

We now assume $\mathcal R>1$ and analyze the stability of the equilibria $E^*=(S^*,I^*)$ by means of the well known 
trace-determinant criterion. Let us consider the 
the Jacobian matrix of the flow of system (\ref{redsys})  at the equilibrium point which, 
taking into account that $\Phi(S^*)=I^*=\Psi(S^*)$, simplifies 
to 
  \[
\mathcal J = \left(\begin{array}{cc}
      -c_{SS}S^*-\dfrac{(\bar{a}r+\bar{\gamma})I^*}{S^*} & \bar{a}r+\bar{\gamma}-(\bar{C}_{SI}+\bar\beta)S^*\\
      (\bar\beta-\bar c_{IS})I^* & -\bar c_{II}I^*
  \end{array}\right).
  \]
This immediately yields $\text{tr} \mathcal F<0$. Furthermore, a direct calculation 
leads to 
\[
 \dfrac{\text{det} \mathcal J}{I^*}  =
\dfrac{\left(c_{SS}\bar c_{II}+(\bar{C}_{SI}+\bar\beta)(\bar\beta-\bar c_{IS})\right)(S^*)^2
+(\bar{a}r+\bar{\gamma})\left(\bar c_{II}I^*-(\bar\beta-\bar c_{IS})S^*\right)}{S^*}
\]
Using again the fact that $I^*=\Psi(S^*)$ if and only if $\bar c_{II}I^*=(\bar\beta-\bar c_{IS})S^*-(\bar \mu+\bar\gamma)$, we have
\[
\dfrac{\text{det}\mathcal J}{I^*}=
\dfrac{\left(k\bar c_{II}+(\bar{C}_{SI}+\bar\beta)(\bar\beta-\bar c_{IS})\right)(S^*)^2
-(\bar{a}r+\bar{\gamma})(\bar \mu+\bar\gamma)}{S^*},
\]
so that $\text{det} \mathcal J>0$ is equivalent to
\begin{equation}\label{c1}
S^*>\sqrt{\dfrac{(\bar{a}r+\bar{\gamma})(\bar \mu+\bar\gamma)}
{k\bar c_{II}+(\bar{C}_{SI}+\bar\beta)(\bar\beta-\bar c_{IS})}},
\end{equation} 
which entails local stability. On the contrary, $E^*$ is unstable if 
\begin{equation}\label{c2}
S^*<\sqrt{\dfrac{(\bar{a}r+\bar{\gamma})(\bar \mu+\bar\gamma)}
{k\bar c_{II}+(\bar{C}_{SI}+\bar\beta)(\bar\beta-\bar c_{IS})}}.
\end{equation} 
The $S$ component of the equilibrium points $E_{+1}^*$ and $E_{+2}^*$ can be explicitly 
calculated from \eqref{roots}.  Direct calculations show that $S_{+1}^*$ fulfills condition \eqref{c2} while condition \eqref{c1} holds for $S_{+2}^*$. \hspace{3cm}
$\blacksquare$\\\newline

\section{Discussion}\label{discu}

We have set up a model aimed to ascertaining the impact of an opportunistic
disease outbreak in a population already affected by a primary disease by assuming that 
both diseases evolve within different time scales. 
For the discussion of results, let us remind the two main aims stated in the introduction. 
On the one hand, we  wanted to know whether the coinfection by a secondary disease produces 
epidemiological scenarios not allowed by the primary disease submodel or not. In the latter case, 
it is of interest to determine if coinfection has any influence on the actual outcome of the model, even if just among those allowed 
by the primary disease submodel.

The answer to the first question is negative, as we have pointed out at the end of section \ref{mod}. 
Thus, the catalog of possible qualitative epidemic behaviors remains unchanged 
by the influence of a secondary disease under the assumptions considered here. 
We can restate this fact by saying that there is neither functional nor dynamical emergence \cite{Auger08a}. 

\begin{figure}[h!]
\begin{center}
\includegraphics[width=12cm]{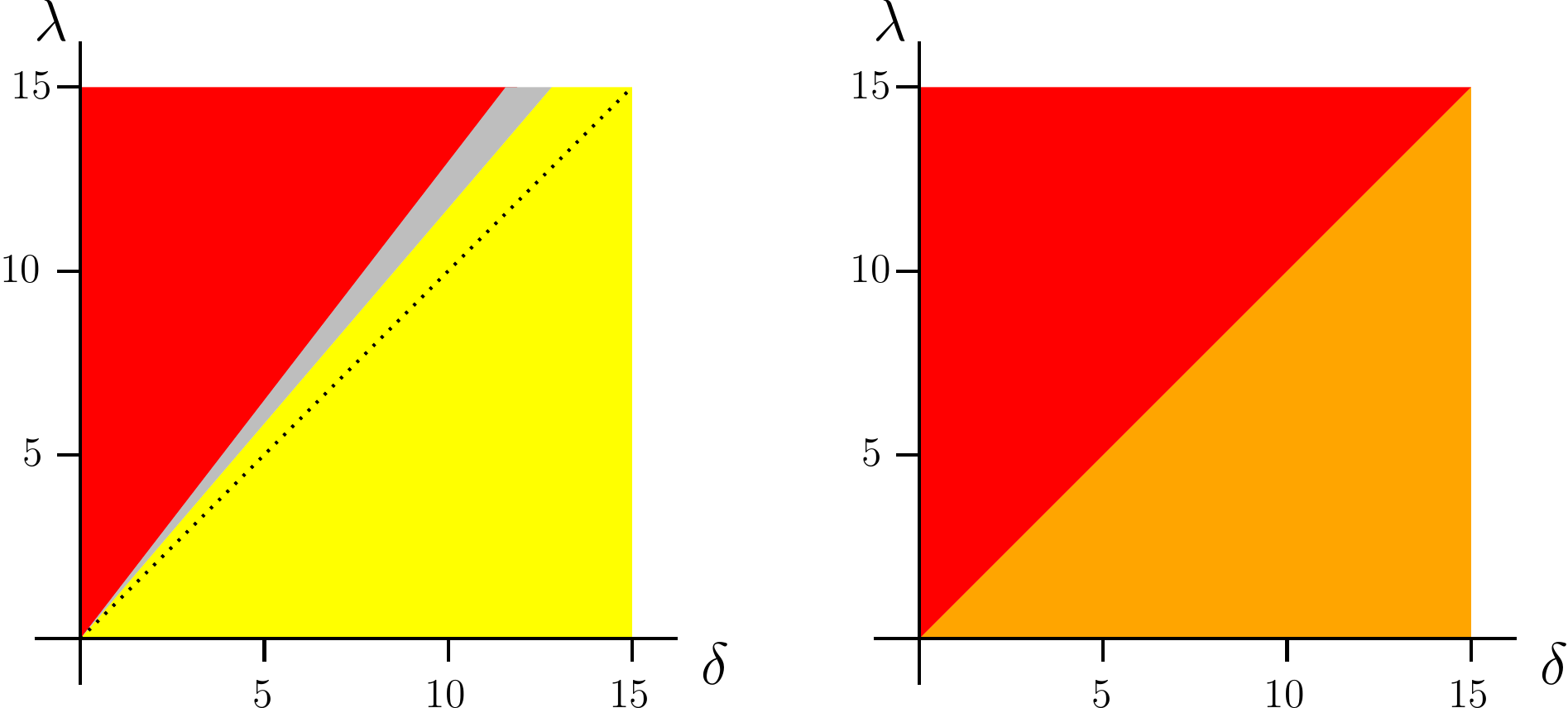}
\caption{Possible epidemic outcomes. 
{Yellow}: disease-free. {Orange}: endemic primary infection. 
{Gray}: disease-free or endemic coinfection, depending on initial values. 
{Red}: endemic coinfection. 
Left figure parameter values:
$a_U=0$.$9$, $a_V=0$.$7$, $r=26$,  $m=12$, $\mu_U=0$.$3$, $\mu_V=0$.$5$, $c_{SS}=3$.$8$, 
$c_{SU}=0$.$5$, $c_{SV}=0$.$5$, $c_{US}=2$.$6$, $c_{UU}=0$.$1$, $c_{UV}=1$, $c_{VS}=0$.$5$,
 $c_{VU}=4 $,  $c_{VV}=4$,  $\beta_U=4 $,  $\beta_V=8 $, $\gamma=0$.$2$. 
The parameter values in the right figure are the same as before but $m=17$, $c_{SS}=2$.$8$, $\beta_V=4 $.}
\end{center} 
\label{picture1}
\end{figure}

Nevertheless, the effect of the the opportunistic disease must be taken into account. 
In Section \ref{anrs} we have found that $\bar A$, $\bar B$, $S_1^*$,  
as well as $\mathcal R$ are key parameters to describe the outcome of the model. 
All of them, but $S_1^*$, depend on 
$\nu^*$, the fraction of coinfected individuals which, in turn, depends on $\lambda/\delta$, 
the ratio of the parameters describing the opportunistic disease dynamics. 
It means that the opportunistic disease can be decisive in the long-term behavior 
of the slow-fast model. 
Therefore, the interest lies on how $\bar A$, $\bar B$ 
and $\mathcal R$ vary with the quotient $\delta/\lambda$. Unfortunately, 
such a dependence is, in general, not simple (just see equations (\ref{fasteq}) and 
(\ref{coeffs})) and we resort to numerical tools to illustrate the effect of varying $\delta/\lambda$ in the outcome of the model. 
 Figure \ref{picture1} displays the different outcomes of the aggregated model as function of $\delta$ and $\lambda$:
In yellow: values of $\delta$ and $\lambda$ leading to the disease-free scenario. 
In orange: values leading to an endemic primary infection outcome.
In gray: conditional coinfection, meaning values leading to either disease-free or endemic coinfection scenario, depending on the the initial amount of susceptible and infected individuals. 
In red: values leading to disease endemic coinfection outcome. 

In the right picture of Figure \ref{picture1}, the epidemiological outcome changes from 
disease-free to endemic coinfection as the ratio $\delta/\lambda$ increases and crosses 
the threshold $\delta/\lambda=1$.  Instead, in the left picture of Figure \ref{picture1}, note that the disease-free region (in the parameters space) overlaps the region $\lambda>\delta$ 
(above the dotted line) where the opportunistic disease would be able to invade if there were primary infected in the population. As the ratio $\delta/\lambda$ increases, the epidemiological outcome changes from disease-free to conditional coinfection and a further increase leads to endemic coinfection. 
The only difference between the parameter values used in each figure is on $c_{SS}$ and $\beta_V$. 
And this fact leads us to another interesting finding: there is a delicate interplay between competition coefficients and infection parameters,
captured by the definition of $\bar A$, $\bar B$, $S_1^*$ and $\mathcal R$, 
which must be taken into account. Although we could not derive general results of such an interdependence, 
we have shown that it must be taken into account.

Summing up, both the irruption of an opportunistic disease and the 
competitive pressure of individuals being in different epidemiological state may affect the 
evolution of the primary disease outbreak. The effect can be determined by means of the parameters  
$\bar A$, $\bar B$ and $\mathcal R$ on $\delta/\lambda$. An it leads us to the second objective 
of this work.
\newline

Related to our second objective, our results point out two different kinds of mechanisms to modulate 
the outcome of the model, each of them feasible within certain ranges of the parameter values. 

On the one hand, having control on parameters $\delta$ and $\lambda$ may allow certain leeway 
to reverse/promote epidemic outbreaks or infection/coinfection scenarios. 
Indeed, once the actual parameter values of the model are known 
one can compute (the equivalent of) Figure \ref{picture1} and get enhanced comprehension on the epidemiological context as well as ascertaining 
the effect on the epidemic outbreak of changes on $\delta$ or $\lambda$. 
In this sense, it is interesting to note that any action or measure taken to modify 
the secondary infection recovery rate $\delta$ or transmission rate $\lambda$ such that 
$\delta/\lambda$ remains constant is completely ineffective. In addition, controlling individuals 
competitive pressure may be relevant for the epidemiological outcome.

On the other hand, the results in propositions \ref{prop3} and \ref{prop6} suggest that acting on the susceptible/infected individuals population size in order to keep the population above/below 
certain threshold allows to have control on epidemic outbreaks. In particular, according to Proposition \ref{prop3}, $\bar A$ is a susceptible population size threshold allowing or not 
the increase of the infection when it is rare. Therefore, introducing/culling (removing) susceptible individuals 
to bring the population above/below this threshold may certainly modify the outcome. 
Besides, under the hypotheses of Proposition \ref{prop6}, we show that when $\mathcal R>1$ whether the disease establishes itself or not depends on the initial amount of susceptible and infected individuals. 
From a mathematical point of view, this scenario is characterized by 
the fact that the disease-free equilibrium $E_1^*$ and an endemic disease (interior) equilibrium $E_{+2}^*$ coexist and are locally asymptotically stable. There is also a saddle node interior equilibrium $E_{+1}^*$. The stable manifold of $E_{+1}^*$ 
separates the basins of attraction of the disease-free and the endemic disease steady states. 
This stable manifold cannot be calculated straightforward, but can be computed using, for instance, the results in \cite{CDRPV14-1} and 
\cite{CDRPV14-2}.

A final comment has to do with the selection of the transmission form of the opportunistic disease. Preliminary calculations show that considering 
DDT instead of FDT leads to equivalent results. This means that even if the nullclines are different, 
the possible outcomes (say the dynamical scenarios) of the corresponding aggregated model are the same.

 \textbf{Acknowledgements.--}
 M. Marv\'a and R. Bravo de la Parra
 are partially supported by Ministerio de Ciencia e Innovaci\'on (Spain),
 projects MTM2011-24321 and MTM2011-25238. E. Venturino is partially supported by
 the project ``Metodi numerici in teoria delle popolazioni''
 of the Dipartimento di Matematica ``Giuseppe Peano''.

\end{document}